\documentclass{article}

\usepackage{latexsym}
\usepackage[centertags]{amsmath}
\usepackage{amsfonts}
\usepackage{amssymb}
\usepackage{amsthm}
\usepackage{newlfont}
\newcommand{\Real}{\mathbb R}
\newcommand{\Complex}{\mathbb C}
\newcommand{\Proj}{\mathbb P}
\newcommand{\Int}{\mathbb Z}
\newcommand{\Nat}{\mathbb N}
\begin{document}

\noindent{\small \bf Ordinary differential equations}\\
\noindent{\small \bf Mechanics of particles and systems}

 \vspace{0.2cm}

\begin{center}
{\LARGE \bf The meromorphic non-integrability of the three-body
problem}\\

 \vspace{0.3cm}

\ {\Large \bf  Tsygvintsev Alexei}

\end{center}

\begin{abstract}
We study the planar three-body problem and prove the absence of a
complete set of complex meromorphic first integrals in a
neighborhood of the Lagrangian solution.
\end{abstract}

\begin{center}
{\bf 1. Introduction}
\end{center}

The three-body problem is a mechanical system which consists of
three mass points $m_1$, $m_2$, $m_3$ which attract each other
according to the Newtonian law [16].

 The practical importance
of this problem arises from its applications to  celestial
mechanics: the bodies which constitute the solar system attract
each other according to  Newton's low, and the stability of this
system on a long period of time is a fundamental question.
Although Sundman [21] gave a power series solution to the
three-body problem in 1913, it was not useful in determining the
growth of the system for long intervals of time. Chazy [3]
proposed in 1922 the first general classification of motion as $t
\rightarrow \infty$. In view of the modern analysis [7], this
stability problem leads to the problem of integrability of a
Hamiltonian system i.e. the existence of a full set of analytic
first integrals in involution. Poincar\'e [18] considered
Hamiltonian functions $H(z,\mu)$ which in addition to
$z_1,\ldots,z_{2n}$ also depended analytically on a parameter
$\mu$ near $\mu=0$. His theorem states that under certain
assumptions about $H(z,0)$, which are in general satisfied, the
Hamiltonian system corresponding to $H(z,\mu)$ can have no
integrals represented as convergent series in $2n+1$ variables
$z_1,\ldots,z_{2n}$ and $\mu$, other than the convergent series in
$H$, $\mu$. Based on this result he proved in 1889 the
non-integrability of the restricted three-body problem [22].
However, this theorem does not assert anything about a fixed
parameter value $\mu$.

 Bruns [2] showed in 1882 that the classical integrals are the only
 independent algebraic integrals of the problem of three bodies.
His theorem has been extended by Painlev\'e [17], who has shown
that every integral of the problem of $n$ bodies which involves
the velocities algebraically (whether the coordinates are involved
algebraically or not) is a combination of the classical integrals.

However, citing [7] `` One may agree with Winter [25] that these
elegant negative results have no importance in dynamics, since
they do not take into account the peculiarities of the behavior of
phase trajectories. As far as first integrals are concerned,
locally, in a neighborhood of a non--singular point, a complete
set of independent integrals always exists. Whether they are
algebraic or transcendent depends explicitly on the choice of
independent variables. Therefore, the problem of the existence of
integrals makes sense only when it is considered in the whole
phase space or in a neighborhood of the invariant set ... ''

 Consider a complex-analytic symplectic  manifold $M$, a
holomorphic Hamiltonian vector field $X_H$ on $M$ and a
non-equilibrium integral curve $\Gamma \subset M$. The nature of
the relationship between the branching of solutions of a system of
variational equations along $\Gamma$ as functions of the complex
time and the non-existence of first integrals of $X_H$ goes back
to the classical works of Kowalewskaya [6]. Ziglin [27] studied
 necessary conditions for an analytic Hamiltonian system with $n
>1$ degrees of freedom to possess $n$ meromorphic independent
first integrals in a sufficiently small neighborhood of the phase
curve $\Gamma$. One can consider the monodromy group $G$ of the
normal variational equations along $\Gamma$. The key idea was that
$n$ independent meromorphic integrals of $X_H$ must induce $n$
independent rational invariants for $G$. Then, in order that
Hamilton's equations have the above first integrals, it is
necessary that for any two non-resonant transformations
$g,g\prime\in G$  $g$ must commute with $g\prime$. Although Ziglin
formulated his result in terms of the monodromy group, it became
quite recently [15,20] that much more could be achieved, under
mild restrictions, by replacing this with the differential Galois
group. Namely, one should check if its identity component, under
Zariski's topology, is abelian.

The collinear three-body problem was proved to be non-integrable
near triple collisions   by Yoshida [26] based on Ziglin's
analysis.

The present paper is devoted to the non-integrability of the
planar three-body problem.

In 1772 Lagrange [8] discovered the particular solution in which
three bodies  form an equilateral triangle and each body describes
a conic.

Moeckel [14] has shown that for a small angular momentum there
exist orbits homoclinic to the Lagrangian elliptical orbits and
heteroclinic between them. Consequently in this case the problem
is not-integrable. Nevertheless, it was observed that for a large
angular momentum and for certain masses of two bodies which are
relatively small compared to the third one, the circular
Lagrangian orbits are stable and, a priori, the system can be
integrable near these solutions. Topan [23] found some examples of
such transcendental integrals in certain configurations of the
restricted three-body problem.

Our approach consists of applying the methods related to [27,15]
to the Lagrangian parabolic orbits.  This means that we will study
the integrability of the problem in a sufficiently small complex
neighborhood of these solutions.

The plan of the paper is follows. In Section 2, following
Whittaker,  we introduce the reductions of the planar three-body
problem from the Hamiltonian system of 6 degrees  of freedom to 3
degrees of freedom. Section 3 is devoted to a parametrization of
the Lagrangian parabolic solution. \\ Section 4 contains the
normal variational equations along this solution. In Section 5 we
study the monodromy group of these equations. In Section 6,
applying the Ziglin's method, we prove that for the three-body
problem there are no two additional meromorphic first integrals in
a connected neighborhood of the Lagrangian parabolic solution
(Theorems 6.2-6.3). Section 7 contains a dynamical interpretation
of above theorems in connection with a theory of splitting and
transverse intersection of asymptotic manifolds.

 \begin{center}
{\bf  2. The reduction of the problem}
\end{center}

Following Whittaker [24] let $(x_1,x_2)$ be the coordinates of
$m_1$, $(x_3,x_4)$ the coordinates of $m_2$, and $(x_5,x_6)$ the
coordinates of $m_3$. Let $y_r=m_k\displaystyle\frac{dx_r}{dt}$,
where $k$ denotes the greatest integer in
$\displaystyle\frac{1}{2}(r+1)$. The equations of  motion are $$
\displaystyle\frac{dx_r}{dt}=\displaystyle\frac{\partial
H_1}{\partial y_r}, \quad
\displaystyle\frac{dy_r}{dt}=-\displaystyle\frac{\partial
H_1}{\partial x_r}, \quad (r=1,2,\dots,6), \leqno (2.1) $$ where
$$
\begin{array}{ll}
H_1=\displaystyle\frac{1}{2m_1}(y_1^2+y_2^2)+\displaystyle\frac{1}{2m_2}(y_3^2+y_4^2)+\displaystyle\frac{1}{2m_3}(y_5^2+y_6^2)-
m_3m_2\{(x_3-x_5)^2+
(x_4-x_6)^2\}^{-1/2}\\-m_3m_1\{(x_5-x_1)^2+(x_6-x_2)^2\}^{-1/2}-
m_1m_2\{(x_1-x_3)^2+ (x_2-x_4)^2\}^{-1/2}.
\end{array}
$$ This is a Hamiltonian system with $6$ degrees of freedom which
admits $4$ first integrals:

\noindent $T_1=H_1 $ -- the energy,\\
$T_2=y_1+y_3+y_5$, $T_3=y_2+y_4+y_6$  -- the components of the impulse of the system,\\
$T_4=y_1x_2+y_3x_4+y_5x_6-x_1y_2-x_3y_4-x_5y_6$  -- the integral of  angular momentum of the system.

The system (2.1) can be transformed to a system with $4$ degrees
of freedom by the following canonical change (Poincar\'e, 1896) $$
x_r=\displaystyle\frac{\partial W_1}{\partial y_r},\quad
g_r=\displaystyle\frac{\partial W_1}{\partial l_r},\quad
(r=1,2,\dots,6), $$ where $$
W_1=y_1l_1+y_2l_2+y_3l_3+y_4l_4+(y_1+y_3+y_5)l_5+(y_2+y_4+y_6)l_6.
\leqno (2.2) $$ Here $(l_1,l_2)$ are the coordinates of $m_1$
relative to axes through $m_3$ parallel to the fixed axes,
$(l_3,l_4)$ are the coordinates of $m_2$ relative to the same
axes, $(l_5,l_6)$ are the coordinates of $m_3$ relative to the
original axes, $(g_1,g_2)$ are the components of impulse of $m_1$,
$(g_3,g_4)$ are the components of impulse of $m_2$, and
$(g_5,g_6)$ are the components of impulse of the system. It can be
shown that in the system of the center of masses the corresponding
equations for $l_5$, $l_6$, $g_5$, $g_6$ disappear from the system
and the reduced system takes the following form $$
\displaystyle\frac{dl_r}{dt}=\displaystyle\frac{\partial
H_2}{\partial g_r},\quad
\displaystyle\frac{dg_r}{dt}=-\displaystyle\frac{
\partial H_2}{\partial l_r}, \quad (r=1,2,3,4), \leqno (2.3)
$$ with the Hamiltonian $$
\begin{array}{ll}
H_2=\displaystyle\frac{M_1}{2}(g_1^2+g_2^2)+\displaystyle\frac{M_2}{2}(g_3^2+g_4^2)+
\displaystyle\frac{1}{m_3}(g_1g_3+g_2g_4)- \displaystyle \frac{
m_3m_2}{\rho_1}- \displaystyle \frac{
m_1m_3}{\rho_2}+\displaystyle\frac{m_1m_2}{\rho_3},
\end{array}
$$ where $$ \rho_1=\sqrt{l_3^2+l_4^2}, \quad
\rho_2=\sqrt{l_1^2+l_2^2}, \quad
\rho_3=\sqrt{(l_1-l_3)^2+(l_2-l_4)^2},$$ are the mutual distances
of the bodies and  $M_1=m_1^{-1}+m_3^{-1}$,
$M_2=m_2^{-1}+m_3^{-1}$.

This system admits two first integrals in involution\\ $K_1=H_2$
-- the energy,\\ $K_2=g_2l_1+g_4l_3+g_6l_5-g_1l_2-g_3l_4-g_5l_6=k$
-- the integral of angular momentum.

Let us suppose that the Hamiltonian system (2.3) possesses a first
integral  $K$ different from  $K_{1,2}$.

\vspace{0.5cm}

\noindent {\bf Definition 2.1} The first integral $K$ of the
system (2.3) is called {\it meromorphic} if it is representable as
a ratio $$ K=\displaystyle \frac{R(l,g)}{Q(l,g)},$$ where $R$, $Q$
are analytic functions of the variables $l_i$, $g_i$, $1 \leq i
\leq 4$.

\vspace{0.5cm}

It can be shown [24] that the system (2.3) possesses an ignorable
coordinate which will make possible a further reduction.

Let us make the following canonical transformation $$
l_r=\displaystyle\frac{\partial W_2}{\partial g_r},\quad
p_r=\displaystyle\frac{\partial W_2}{\partial q_r},\quad
(r=1,2,3,4), \leqno (2.4) $$ where $$
W_2=g_1q_1\mathrm{cos}q_4+g_2q_1
\mathrm{sin}q_4+g_3(q_2\mathrm{cos}q_4-
q_3\mathrm{sin}q_4)+g_4(q_2\mathrm{sin}q_4+q_3\mathrm{cos}q_4). $$
Here $q_1$ is the distance $m_3m_1$; $q_2$ and $q_3$ are the
projections of $m_2m_3$ on, and perpendicular to $m_1m_3$; $p_1$
is the component of momentum of $m_1$ along $m_3m_1$; $p_2$ and
$p_3$ are the components of momentum of $m_2$ parallel and
perpendicular to $m_3m_1$.

One can write the new equations as follows
$$
\displaystyle\frac{dq_r}{dt}=\displaystyle\frac{\partial H}{\partial p_r},\quad \displaystyle\frac{dp_r}{dt}=
-\displaystyle\frac{
\partial H}{\partial q_r}, \quad (r=1,2,3), \leqno (2.5)
$$ and $$ \displaystyle\frac{dq_4}{dt}=\displaystyle\frac{\partial
H}{\partial p_4},\quad \displaystyle\frac{dp_4}{dt}=0, \leqno
(2.5.a) $$ with the Hamiltonian $$
\begin{array}{ll}
H=\displaystyle \frac{M_1}{2}\left\{p_1^2+\displaystyle
\frac{1}{q^2_1}P^2\right\}+\displaystyle \frac{M_2}{2}(p_2^2+
p_3^2)+\displaystyle \frac{1}{m_3}\left\{p_1p_2-\displaystyle
\frac{p_3}{q_1}P\right\} -\displaystyle \frac{ m_1m_3}{r_1}-
\displaystyle \frac{m_3m_2}{r_2}- \displaystyle \frac{
m_1m_2}{r_3}, \\ P=p_3q_2-p_2q_3-p_4,
\end{array}
$$  where $$ r_1=q_1, \quad r_2=\sqrt{q^2_2+q^2_3}, \quad
r_3=\sqrt{ (q_1-q_2)^2+q^2_3},$$ are the mutual distances of the
bodies.

Since $p_4=k=const$ the system (2.5) is a closed Hamiltonian
system with $3$ degrees of freedom.  If this system is integrated
then $q_4$ can be found by a quadrature from (2.5.a).

\vspace{0.5cm}

\noindent {\bf Proposition 2.2} {\it If the Hamiltonian system
(2.3) admits the full set of functionally independent meromorphic
first integrals in involution $\{ K_1,K_2,K_3,K_4\}$  then the
system (2.5) possesses two functionally independent additional
first integrals $\{H_1,H_2\}$ which are  meromorphic functions of
the variables $q_i$, $p_i$, $1\leq i \leq 3$.}

This is the obvious consequence of  the canonical change (2.4).

\begin{center}
{\bf 3. A parametrization of the parabolic Lagrangian solution}
\end{center}

The equations (2.1) admit an exact  solution discovered by
Lagrange [8] in which the triangle formed by the three bodies is
equilateral and the trajectories of the bodies are similar conics
with one focus at the common barycenter.
 For the reduced form (2.5) the equality of the mutual distances gives
$$ q_1=q,\quad q_2=\displaystyle\frac{q}{2},\quad
q_3=\displaystyle\frac{\sqrt{3}q}{2}, \leqno (3.1) $$ where
$q=q(t)$ is an unknown function. Substituting (3.1) into (2.5) one
can show that $$ p_1=p,\quad p_2=Ap+\displaystyle\frac{B}{q},\quad
p_3=Cp+\displaystyle\frac{D}{q}, \leqno (3.2) $$ with $p=p(t)$
unknown and  $A$, $B$, $C$, $D$ are the following constants $$
\begin{array}{ll} A=\displaystyle \frac{m_2(m_3-m_1)}{m_1S_3},\quad
B=-\displaystyle \frac{\sqrt{3}kS_1m_2m_3}{S_2S_3},\quad
C=\displaystyle \frac{\sqrt{3}m_2(m_1+m_3)}{m_1S_3},\\
D=-\displaystyle \frac{km_2(S_2+m_1m_2-m_3^2)}{S_2S_3},
\end{array} $$
 where
 $$ S_1=m_1+m_2+m_3, \quad
S_2=m_1m_2+m_2m_3+m_3m_1, \quad S_3=m_2+2m_3.$$

Substituting (3.1), (3.2) into  the integral of energy $H=h=const$
we obtain the following relation between $q$ and $p$ $$
ap^2+\displaystyle\frac{bp}{q}+\displaystyle\frac{c}{q}+\displaystyle\frac{d}{q^2}=h,
\leqno (3.3) $$ where $$
a=\displaystyle\frac{2S_1S_2}{m_1^2S_3^2},\quad
b=-\displaystyle\frac{2\sqrt{3}km_2S_1}{m_1S_3^2},\quad c=-S_2,
\quad d=\displaystyle\frac{2k^2S_1(m_2^2+m_2m_3+m_3^2)}{S_3^2S_2}.
$$ Moreover, from  (2.5) we have $$
\displaystyle\frac{dq}{dt}=\left(M_1+\displaystyle\frac{A}{m_3}\right)p+\displaystyle\frac{B}{m_3q}
 \leqno (3.4)
$$ The equations (3.1), (3.2), (3.3), (3.4) define all Lagrangian
particular solutions and contain two free parameters: $k$ and $h$.

Consider the case of zero energy $h=0$ and $k\neq0$. Then  there
exists a parabolic particular solution in the sense that the limit
velocity goes to zero when the bodies approach infinity and each
body describes a parabola.

Putting $w=pq$  one can find  by using of (3.3) $q$, $p$ as the
functions of $w$ $$ q=P(w),\quad p=\displaystyle\frac{w}{P(w)},
\leqno (3.5) $$
 where $P(w)=-(aw^2+bw+d)/c$.

Let $M={\Complex^6}$ be the complexified phase space of the system
(2.5). Then (3.5), (3.1), (3.2) define a parametrized parabolic
integral curve $\Gamma \in M$ with the parameter $w\in {\Complex
\Proj^1}$.

\begin{center}
{\bf 4. The normal  variational equations}
\end{center}

Let $z=(q_1,q_2,q_3,p_1,p_2,p_3)$, $z\in M$. One can obtain the
variational equations of the system (2.5)  along the integral
curve $\Gamma$ $$
\displaystyle\frac{d\zeta}{dt}=JH_{zz}(\Gamma)\zeta, \quad \zeta
\in T_{\Gamma}M,  \leqno (4.1) $$ where  $H_{zz}$ is the Hessian
matrix of Hamiltonian $H$ at $\Gamma$ and $J$ is the $6\times 6$
matrix $$ J
=
\left (\begin {array}{cc} 0&E\\\noalign{\medskip}-E&0\end {array}
\right ),
 $$ where $E$ is the identity $3 \times 3$  matrix .

These equations admit the linear first integral
$F=(\zeta,H_z(\Gamma))$, where $H_z=grad(H)$ and can be reduced on
the normal $5$-dimensional bundle $ G=T_{\Gamma}M/T\Gamma$ of
$\Gamma$ . After the restriction of (4.1) on the surface $F=0$ we
obtain {\it normal variational equations} (NVE) [27] which are the
system of $4$ equations $$ \displaystyle\frac{d\eta}{dt}=\tilde
A(\Gamma)\eta, \quad \eta\in{\Complex^4}, \leqno (4.2) $$ where
$\tilde A$ is a  $4 \times 4$ matrix depending on $\Gamma$.

We can obtain NVE in the following natural way applying
Whittaker's procedure [24] of reducing the order of the
Hamiltonian system (2.5).

Fixing the level of energy $h=0$ one can  find $p_1$ as a function
of the other variables from  the equation $H(q,p)=0$ which takes
the following form  $$ a_1p_1^2+b_1p_1+c_1=0, $$ where $a_1$,
$b_1$, $c_1$ are  known functions depending on $p_2$, $p_3$,
$q_1$, $q_2$, $q_3$.

Solving this  equation we get two solutions for $p_1$ $$
p_1=\displaystyle\frac{-b_1+\sqrt{\Delta}}{2a_1}=K_+ \quad
\mathrm{and} \quad
p_1=\displaystyle\frac{-b_1-\sqrt{\Delta}}{2a_1}=K_- ,$$ where
$\Delta=b_1^2-4a_1c_1$.

By substituting  the Lagrangian solution given by (3.1), (3.2),
(3.5) in these relations we choose the root $p_1=K_-$ as
corresponding to this solution.

The functions $q_r(t)$, $p_r(t)$, $r=2,3$ satisfy the canonical
equations $$
\displaystyle\frac{dq_r}{dq_1}=\displaystyle\frac{\partial
K}{\partial p_r},\quad
\displaystyle\frac{dp_r}{dq_1}=-\displaystyle\frac{
\partial K}{\partial q_r}, \quad (r=2,3), \leqno (4.3)
$$ where $K=-K_-$ and $q_1$ is taken as the new time.

The system (4.3) is a nonautonomous Hamiltonian system with $2$
degrees of freedom which has the same integral curve $\Gamma$.
Notice that $K$ is not more a first integral.

It is useful to pass now to the new time $(q_1=q)\rightarrow w$.
From the formulas (3.3), (3.5) we have $$
q=\displaystyle\frac{aw^2+bw+d}{c},\quad
dq=-\displaystyle\frac{2aw+b}{c}dw. \leqno (4.4) $$ The resulting
NVE (4.2) are obtained as the variational equations of the system
(4.3) near the integral curve $\Gamma$ and after the substitution
(4.4) take the form $$ \displaystyle\frac{d\eta}{dw}=\tilde
A(\Gamma)\eta, \quad \eta\in{\Complex^4}, \leqno (4.5) $$ where
$\tilde A$ is a $4\times 4$ matrix whose elements are rational
functions of $w$.

We can represent $\tilde A$ in the following block form
$$
\tilde A
=
\left (\begin {array}{cc}
{M_{{3}}}^{T}&M_{{2}}\\\noalign{\medskip}-M_ {{1}}&-M_{{3}}\end
{array}\right ), $$ where $M_1$, $M_2$, $M_3$ are $2\times 2$
matrices and $M_3^T$ means the transposition of $M_3$.

The matrix $M_1$ is symmetric and has the following form
$$
M_1
=
\displaystyle\frac{1}{S_1^3L^2Z^2}
\left(
\begin{array}{ll}
n_{11}&n_{12} \\
n_{12}&n_{22}
\end{array}
\right),
$$
where $L(w)$ is the linear polynomial
$$
L=l_1w+l_2,
$$
and $l_1=2S_2$,\quad $l_2=-\sqrt{3}m_1m_2k$.

$Z(w)$ is the following quadratic polynomial
 $$ Z=z_1w^2+z_2w+z_3,
$$
 where $z_1=S_2^2,\quad z_2=-\sqrt{3}m_1m_2kS_2,\quad
z_3=k^2m_1^2(m_2^2+m_2m_3+m_3^2).$

The coefficients $n_{ij}$ have the form $$
n_{11}=A_1w^2+A_2w+A_3,\quad n_{12}=A_4w^2+A_5w+A_6,\quad
n_{22}=A_7w^2+A_8w+A_9, $$ where $A_i$ are  constants depending on
the masses $m_1$, $m_2$, $m_3$ and  $k$.

The matrix $M_2$ has the following expression
$$
M_2
=
\displaystyle\frac{4S_1Z}{S_2S_3^3m_1^4m_2m_3} \left(
\begin{array}{ll}
1&0 \\
0&1
\end{array}
\right).
$$
For the matrix $M_3$ we have
$$
M_3
=
\displaystyle\frac{1}{m_1S_1LZ}
\left(
\begin{array}{ll}
m_{11}&m_{12} \\
m_{21}&m_{22}
\end{array}
\right),
$$
where
$$
\begin{array}{ll}
m_{11}=B_1w^2+B_2w+B_3,\quad m_{12}=B_4w^2+B_5w+B_6,\quad m_{21}=B_7w^2+B_8w+B_9,\\
m_{22}=B_{10}w^2+B_{11}w+B_{12},
\end{array}
$$ and $B_j$ are  constants depending on $m_1$, $m_2$, $m_3$ and
$k$.

The system (4.5) has four singular points $w_1$, $w_2$, $w_3$,
$w_4$ in the complex plane: $$ w_1=\infty, $$ -- the infinity. $$
w_2=\displaystyle\frac{\sqrt{3}m_1m_2k}{2S_2}, $$ -- the root of
$L=0$. $$ w_3=\displaystyle \frac {( \sqrt{3}m_2 + i S_3)km_1}{2
S_2} ,\quad w_4=\displaystyle \frac{(\displaystyle  \sqrt{3}m_2 -
\displaystyle  i S_3)k m_1}{2S_2}, \leqno (4.6) $$ -- the
corresponding roots of the quadratic equation $Z=0$ where
$i^2=-1$.

Notice that  the expressions for $w_{2,3,4}$ have a rational form
on the masses.

The singularities $w_i$, $1 \leq i \leq 4$  have a clear
mechanical sense: $w_1$ corresponds to the motion of the bodies at
infinity, $w_2$ defines the moment of the  maximal approach.

It is easy to see from (4.6) that if the angular momentum constant
$k=0$, then $w_2=w_3=w_4=0$ and we have a triple collision of the
bodies at the moment of time $w=0$. If $k\neq 0$ then by the lemma
of Sundman there are no triple collisions in the real phase space
and $w_{3,4}$ become complex.

Since the expression for $p$ given in (3.5) becomes infinity when
$w\rightarrow w_{3,4}$, formally, we can consider $w_3$ and $w_4$
as corresponding to the ``complex'' collisions which tend to $w=0$
as $k\rightarrow 0$.

It was noted by Schaefke [19] that the equations (4.5) can be
reduced to  fuchsian form.

In order to do it, consider  the linear change of variables
$$\eta=Cx, \leqno (4.7) $$
where
$\eta=(\eta_1,\eta_2,\eta_3,\eta_4)^T$, $x=(x_1,x_2,x_3,x_4)^T$
and $C=diag(LZ,LZ,1,1)$.

In new variables the system (4.5) takes the following form $$
\displaystyle\frac{dx}{dw}=\left(\displaystyle\frac{A(k)}{w-w_2}+\displaystyle\frac{B(k)}{w-w_3}+
\displaystyle\frac{C(k)}{w-w_4}\right)x, \quad  x\in{\Complex^4},
\leqno (4.8) $$ where $A(k),B(k),C(k)$ are known constant $4\times
4$ matrices depending on $m_1,m_2,m_3$ and $k$.

Under the assumption $k\neq 0$ we can exclude the parameter $k$
from the system (4.8) by using the change of time $w=kt$. As a
result, one obtains $$ \displaystyle\frac{dx}{dt}=\left(
\displaystyle\frac{A}{t-t_0}+\displaystyle\frac{B}{t-t_1}+\displaystyle\frac{C}{t-t_2}\right)x,
\leqno (4.9) $$

where $$ t_0=\displaystyle\frac{\sqrt{3}m_1m_2}{2S_2},\quad
t_1=\displaystyle\frac{m_1(\sqrt{3}m_2+iS_3)}{2S_2}, \quad
t_2=\displaystyle\frac{m_1(\sqrt{3}m_2-iS_3)}{2S_2}. $$ and $$
A=\displaystyle\frac{\tilde M(t_0)}{(t_0-t_1)(t_0-t_2)},\quad
B=\displaystyle\frac{\tilde M(t_1)}{(t_1-t_2)(t_1-t_0)},\quad
C=\displaystyle\frac{\tilde M(t_2)}{(t_2-t_1)(t_2-t_0)}. $$

Here, $\tilde M(w)$ is the following matrix $$ \tilde M(w)= \left(
{\begin{array}{cc} L\,Z\,{M_{3}}^{T} - {\displaystyle \frac {{
\partial L}\,Z}{{ \partial w }}E}  & {M_{2}} \\
 - L^{2}\,Z^{2}\,{M_{1}} &  - L\,Z\,{M_{3}}
\end{array}}
 \right), $$
where one should put $k=1$.

The system (4.9) is defined on a connected Riemann surface
$X={\Complex \Proj^1 / \{t_0,t_1,t_2,\infty \}}$.

It turns out that the matrix $A$ is real and the matrices
$B=R+iJ$, $C=R-iJ$ are complex conjugate being $R$ and  $J$ real
matrices. It will simplify matters further if we choose the units
of masses as follows $$m_1=\alpha,\quad  m_2=\beta, \quad
m_3=1,\quad 0< \alpha\leq\beta\leq 1 .$$ In Appendix A we write
the expressions for $A$, $R$, $J$ with help of MAPLE.

\begin{center}
{\bf  5. The monodromy group of the system (4.9)}
\end{center}

Let $\Sigma(t)$ be a solution of the matrix equation (4.9) $$
\displaystyle\frac{d}{dt}\Sigma=\left(
\displaystyle\frac{A}{t-t_0}+\displaystyle\frac{B}{t-t_1}+\displaystyle\frac{C}{t-t_2}\right)\Sigma,
\leqno (5.1) $$ with the initial condition $\Sigma(\tau)=I$, $\tau
\in X$ where $I$ is the unit $4\times 4$ matrix.

It can be continued along a closed path $\gamma$ with end points
at $\tau$. We obtain the function $\tilde \Sigma(t)$ which also
satisfies (5.1). From linearity of (5.1) it follows that there
exists a complex $4\times 4$ matrix $T_{\gamma}$ such that $\tilde
\Sigma(t)=\Sigma(t) T_{\gamma}$. The set of matrices
$G=\{T_{\gamma}\}$ corresponding to all closed curves in $X$ is a
group. This group is called the {\it monodromy group} of the
linear system (4.9). Let $T_i$ be the elements of $G$
corresponding to circuits around the singular points $t=t_i$,
$i=0,1,2$. Then the monodromy group $G$ is formed by $T_0$, $T_1$,
$T_2$. Denote by $T_{\infty}\in G$ the element corresponding to a
circuit around the point $t=\infty$.

\vspace{0.5cm}

\noindent {\bf Lemma 5.1} {\it The following assertions about the
monodromy group $G$ hold
\\

\noindent a)  $T_0=I$ --  is the unit matrix and $$
T_1T_2=T^{-1}_{\infty}. \leqno (5.2) $$

\vspace{0.5cm}

\noindent b) There exist two non-singular matrices $U$, $V$ such
that $$ U^{-1}T_1U=V^{-1}T_2V=\left(
\begin{array}{cccc}
1 & 1 & 0 & 0 \\ 0 & 1 & 0 & 0 \\ 0 & 0 & 1 & 1
\\ 0 & 0 & 0 & 1
\end{array}
\right).$$

\vspace{0.5cm}

\noindent c) The matrix $T_{\infty}$ has the following eigenvalues
$$\mathrm{Spectr}(T_{\infty})=\left\{ e^{2\pi i\lambda_1},\quad
e^{2\pi i\lambda_2},\quad e^{-2\pi i\lambda_1},\quad e^{-2\pi
i\lambda_2}\right\}, \leqno (5.3)$$ where $$
\lambda_1=\displaystyle\frac{3}{2}+\displaystyle\frac{1}{2}\sqrt{13+\sqrt{\theta}},\quad
\lambda_2=\displaystyle\frac{3}{2}+\displaystyle\frac{1}{2}\sqrt{13-\sqrt{\theta}},
\leqno (5.4) $$ and $$ \theta=144\left(1-\displaystyle
3\frac{S_2}{S_1^2}\right),\quad S_1=\alpha+\beta+1,\quad
S_2=\alpha\beta+\alpha+\beta.$$ Moreover, $$ \mathrm{Spectr}
(T_{\infty} ) \neq \{1,1,1,1\}. $$}

  \vspace{0.5cm}

 {\it Proof.}  a) The matrix $A$ has the eigenvalues
$\{-1,-1,0,0\}$. Following the general theory of the linear
differential equations let us write the general solution of the
system (4.9) near the singular point $t=t_0$ as follows $$
x(t)=c_1X_1(t)+c_2X_2(t)+c_3X_3(t)+c_4X_4(t),$$ where
$c_{1,\ldots,4}\in {\Complex}$ are arbitrary constants and $$
\begin{array}{llcc}
X_1(t)=\displaystyle
\frac{a_{-1}}{t-t_0}+a_0+a_1(t-t_0)+\cdots,\quad
X_2(t)=\displaystyle \frac{b_{-1}}{t-t_0}+b_0+b_1(t-t_0)+\cdots,\\
X_3(t)=c_0+c_1(t-t_0)+\cdots,\quad \quad \quad \quad \quad
X_4(t)=d_0+d_1(t-t_0)+\cdots,
\end{array} \leqno (5.5) $$ where $a_i$, $b_i$, $c_i$, $d_i \in
{\Complex^4} $ are some constant vectors.

By substituting (5.5) in (4.9) one can find  $a_i$, $b_i$, $c_i$,
$d_i$ and show that the vectors $X_1(t)$, $X_2(t)$, $X_3(t)$,
$X_4(t)$ are functionally independent and meromorphic in a small
neighborhood of the point $t=t_0$. This implies that the element
$T_0$ of the monodromy group $G$ corresponding to a circuit around
$t_0$ is the unit matrix. Obviously we should have
$T_0T_1T_2=T^{-1}_{\infty}$. From this fact the relation (5.2)
follows.

\vspace{0.5cm}

\noindent b) The matrices $B$, $C$ have the  same eigenvalues
$\{-2,-1,0,1\}$. It can be shown by a straightforward calculation
that near the singular point $t=t_1$ the general solution of the
system (4.9) can be represented as $$
x(t)=c_1Y_1(t)+c_2Y_2(t)+c_3Y_3(t)+c_4Y_4(t),$$ where
$c_{1,\ldots,4}\in {\Complex}$ are arbitrary constants and $$
\begin{array}{lllcc} Y_1(t)=e_1(t-t_1)+e_2(t-t_1)^2+\cdots, \quad
Y_2=f_0+f_1(t-t_1)+\cdots+C_1 \mathrm{ln}(t-t_1)Y_1(t), \\
Y_3(t)=\displaystyle \frac{g_{-1}}{t-t_1}+g_0+g_1(t-t_1)+\cdots,\\
Y_4(t)=\displaystyle \frac{h_{-2}}{(t-t_1)^2}+\displaystyle
\frac{h_{-1}}{t-t_1}+\cdots+C_2\mathrm{ln}(t-t_1)(f_0+f_1(t-t_1)+\cdots)+C_3\mathrm{ln}(t-t_1)Y_1(t),
\end{array}
$$ where $e_i$, $f_i$, $g_i$, $h_i \in {\Complex^4}$ are some
constant vectors and $C_1$, $C_2$, $C_3$ are parameters depending
on the masses $\alpha$, $\beta$.

  For $C_1$, $C_2$ one can find $$
 C_1=\displaystyle \frac{9}{4}\frac{\beta \alpha^3( \beta+2)^2(\alpha \beta+\alpha+
\beta)}{(\alpha+\beta+1)^3 },\quad C_2=iC_1.$$

The matrix $\Sigma(t)=(Y_1,Y_2,Y_3,Y_4)$ represents the solution
of the system (5.1) in a small neighborhood of the point $t=t_1$.
After going around of $t_1$ we get $\tilde \Sigma(t) =\Sigma(t) M$
where $$ M=\left(
\begin{array}{cccc}
1 & 2\pi i C_1 & 0 & 2\pi i C_3 \\ 0 & 1 & 0 & 0 \\ 0 & 0 & 1 &
2\pi iC_2
\\ 0 & 0 & 0 & 1
\end{array}
\right). $$

Since $C_1 \neq 0$, $C_2\neq 0$ for $\alpha>0$, $\beta>0$, there
exists a non-singular matrix $T$ such that $$ T^{-1}MT=\left(
\begin{array}{cccc}
1 & 1 & 0 & 0 \\ 0 & 1 & 0 & 0 \\ 0 & 0 & 1 & 1
\\ 0 & 0 & 0 & 1
\end{array}
\right), \leqno (5.6) $$ which is the Jordan form of $M$.

The matrix $T_1$ is similar to $M$ and therefore has the same
Jordan form (5.6). Repeating the analogous  arguments for the
matrix $T_2$ we  deduce that the same assertion holds for the
monodromy matrix  $T_2$. Notice that the existence of logarithmic
branching near some Lagrangian solutions in three body problem was
first observed by H. Block (1909) and J.F. Chazy (1918) ( see for
instance [1]).

\vspace{0.5cm}

\noindent c) Consider the matrix $A_{\infty}=-(A+B+C)$. Then there
exists (see for example [4]) a non-singular matrix $W$ such that
$$ T_{\infty}=W^{-1}e^{2\pi iA_{\infty}}W. \leqno (5.7) $$

Appendix A contains the expressions for the elements of the matrix
$A_{\infty}$. One can calculate its eigenvalues $$
\mathrm{Spectr}(A_{\infty})=\{
\lambda_1,\lambda_2,3-\lambda_1,3-\lambda_2\}, $$ where
$\lambda_{1,2}$ are given in (5.3).

One can easy check that $$ 0\leq \sqrt{\theta} < 12, \leqno (5.8)
$$ for all  $\alpha >0$, $\beta > 0$.

With the help of (5.7) we obtain for the eigenvalues of the matrix
$T_{\infty}$ the expression (5.3).

Let us suppose now that $\mathrm{
Spectr}(T_{\infty})=\{1,1,1,1\}$. Then according to (5.4) we
obtain
 $$ \sqrt{13+\sqrt{\theta}}=n_1, \quad \sqrt{13-\sqrt{\theta}}=n_2, \quad n_1,n_2
 \in{\Int
Z}.\leqno (5.9) $$ Hence, in view of (5.8), the number $r
=\sqrt{\theta}$ is an integer $0\leq r \leq 11$. The simple
calculation shows  that for these $r$ the relations (5.9) are not
fulfilled. This implies that $$ \mathrm{Spectr} (T_{\infty} ) \neq
\{1,1,1,1\}. $$ The proof of Lemma 5.1 is completed. \quad \quad
$\Box$

\newpage

\begin{center}

{\bf 6. Nonexistence of additional meromorphic first integrals}

\end{center}

We call the planar three-body problem (2.1)  {\it meromorphically}
integrable near the Lagrangian parabolic solution $\Gamma$,
defined in Section 3, if the corresponding Hamiltonian system
(2.3) possesses a complete set of complex meromorphic first
integrals (see Definition 2.1)  in involution in a connected
neighborhood of $\Gamma$. Recall that equations (2.3) describe the
motion of bodies in the system of the center of masses.

From Proposition 2.2 it follows that in this case the system (2.5)
admits two additional first integrals which are meromorphic and
functionally independent in the same  neighborhood.

\vspace{0.5cm}

\noindent {\bf Theorem 6.1} {\it For $k\neq 0$ for the Hamiltonian
system (2.5) there are no two functionally independent additional
first integrals, meromorphic in a connected neighborhood of the
Lagrangian parabolic solution $\Gamma$.}

\vspace{0.5cm}

{\it Proof.} Suppose that the Hamiltonian system (2.5) admits two
functionally independent first integrals $H_1$, $H_2$, meromorphic
in a connected neighborhood of the Lagrangian parabolic solution
$\Gamma$ and functionally independent together with $H$. According
to Ziglin [27] in this case the NVE (4.5) have two functionally
independent meromorphic integrals $F_1$, $F_2$ which are
single-valued in a complex neighborhood of the Riemann surface
$\Gamma={\Complex \Proj^1}/\{t_0,t_1,t_2,\infty\}$. The linear
system (4.9) was obtained from (4.5) by the linear change of
variables (4.7) and the change  of the time $w=kt$, $k\neq 0$.
Therefore, it possesses two functionally independent meromorphic
integrals $I_1$, $I_2$. From this fact the following lemma is
deduced

\vspace{0.5cm}

\noindent{\bf Lemma 6.2 (Ziglin [27])} {\it The monodromy group
$G$ of the system (4.9) has two rational, functionally independent
 invariants $J_1$, $J_2$.}

\vspace{0.5cm}

 In appropriate coordinates, according to b) of Lemma 5.1, the monodromy transformation $T_1$
can be written as follows $$T_1= \left(
\begin{array}{cccc}
1 & 1 & 0 & 0 \\ 0 & 1 & 0 & 0 \\ 0 & 0 & 1 & 1
\\ 0 & 0 & 0 & 1
\end{array}
\right) =I+D,$$ where $I$ is the unit matrix and $$D=\left(
\begin{array}{cccc}
0& 1 & 0 & 0 \\ 0 & 0 & 0 & 0 \\ 0 & 0 & 0 & 1
\\ 0 & 0 & 0 & 0
\end{array}
\right ). \leqno (6.1) $$

For the monodromy matrix $T_2$ one writes $$ T_2=I+R,$$ where
$$R=\tilde V D \tilde V^{-1}=\left(
\begin{array}{cccc}
a_1 & a_2 & a_3 & a_4 \\ b_1 & b_2 & b_3 & b_4 \\ c_1 & c_2 & c_3
& c_4 \\ d_1 & d_2 & d_3 & d_4
\end{array}
\right),\leqno (6.2)$$ with some unknowns $a_i$, $b_i$, $c_i$,
$d_i\in {\Complex}$ and a nonsingular matrix $\tilde V$.

Let us input the following linear differential operators $$
\delta=x_2\displaystyle\frac{\partial}{\partial
x_1}+x_4\displaystyle\frac{\partial}{\partial x_3},$$ and $$
\Delta=
\left(\sum_{i=1}^4a_ix_i\right)\displaystyle\frac{\partial}{\partial
x_1}+\left(\sum_{i=1}^4b_ix_i\right)\displaystyle\frac{\partial}{\partial
x_2}+\left(\sum_{i=1}^4c_ix_i\right)\displaystyle\frac{\partial}{\partial
x_3}+\left(\sum_{i=1}^4d_ix_i\right)\displaystyle\frac{\partial}{\partial
x_4}. $$

\vspace{0.5cm}

\noindent{\bf Lemma 6.3} {\it Let $J$ be a rational invariant of
the monodromy group $G$, then the following relations hold $$
\delta J=0, \quad \Delta J =0.$$}

\vspace{0.5cm}

{\it Proof.} For an arbitrary $n\in \Nat$ we have $T_1^n=I+nD$,
hence $J\left(T_1^n x\right)=J\left( x+nD x\right).$ Expanding the
last expression in  Taylor series we obtain  $$ J\left(T_1^n
x\right)=J(x)+n\delta J(x)+\sum\limits_{i=2}^{\infty} n^i r_i(x),
\leqno (6.3)$$ where $r_i(x)$ are some rational functions.

In view of $J\left(T_1^n x\right)=J(x)$ and the fact that $J(x)$
is a rational function on $x$,  the second term of (6.3) gives
$\delta J=0$. The relation $\Delta J=0$ is deduced by analogy from
the identity $J\left(T_2x\right)=J(x)$. \quad \quad $\Box$

{\it Case (1).} Assume that  invariants $J_1$, $J_2$ depend on
$x_2$, $x_4$ only. By Lemma 6.3 we have $$ \Delta J_1=0, \quad
\Delta J_2=0. \leqno (6.4)$$

It can be verified that  the equations (6.4) imply  the conditions
$b_i=0$, $d_i=0$, $1\leq i \leq 4$. Accordingly, the matrix $R$
may be written  $$ R=\left(
\begin{array}{cccc}
a_1 & a_2 & a_3 & a_4 \\ 0 & 0 & 0 & 0 \\ c_1 & c_2 & c_3 & c_4
\\ 0 & 0 & 0 & 0
\end{array}
\right). \leqno (6.5)$$

One can find the characteristic polynomial
$P(\lambda)=det(R-\lambda I)$ of  $R$
$$P(\lambda)=\lambda^4-(a_1+c_3)\lambda^3+(a_1c_3-c_1a_3)\lambda^2.
\leqno (6.6)$$

In view of (6.1), (6.2) all eigenvalues of the matrix $R$ are
equal to $0$, thus, with help of (6.6) we get $$ a_1+c_3=0, \quad
a_1c_3=c_1a_3. \leqno (6.7) $$

The matrix $T_1T_2$ takes the following form $$ T_1T_2=\left(
\begin{array}{cccc}
a_1+1 & a_2+1 & a_3 & a_4 \\ 0 & 1 & 0 & 0 \\ c_1 & c_2 & c_3+1 &
c_4+1
\\ 0 & 0 & 0 & 1
\end{array}
\right), $$ and $$ \mathrm{Spectr}(T_1T_2)=\{1,1,s+f,s-f\},$$
where $$ s=1+\displaystyle \frac{a_1+c_3}{2}, \quad
f=\displaystyle \frac{ \sqrt { a_1^2+ c_3^2+4 c_1a_3-2a_1
c_3}}{2}. \leqno (6.8) $$

The straightforward  calculation by using (6.7) and (6.8) shows
that the eigenvalues of the matrix $T_1T_2$ are equal to
$\{1,1,1,1\} $. According to (5.2) these must be the eigenvalues
of the matrix $T_{\infty}$ which is in contradiction to c) of
Lemma 5.1.

{\it Case(2).} Assume that even one from the  invariants $J_1$,
$J_2$ depends on $x_1$ or $x_3$. Let, for example $$
\displaystyle\frac{\partial J_1}{\partial x_1}\neq 0. \leqno (6.9)
$$

It is useful to consider two additional linear operators $
\delta_1=[\delta,\Delta]$ and $
\delta_2=-\displaystyle\frac{1}{2}[\delta,\delta_1].$

One has $$ \delta_1=f_1\displaystyle\frac{\partial}{\partial
x_1}+f_2\displaystyle\frac{\partial}{\partial x_2}+
f_3\displaystyle\frac{\partial}{\partial
x_3}+f_4\displaystyle\frac{\partial}{\partial x_4}, \quad
\delta_2=(b_1x_2+b_3x_4)\displaystyle\frac{\partial}{\partial
x_1}+ (d_1x_2+d_3x_4)\displaystyle\frac{\partial}{\partial x_3},
$$ where $$
\begin{array}{llcc}
f_1=-b_1x_1+(a_1-b_2)x_2-b_3x_3+(a_3-b_4)x_4, &
f_2=b_1x_2+b_3x_4,\\ f_3=-d_1x_1+(c_1-d_2)x_2-d_3x_3+(c_3-d_4)x_4,
& f_4=d_1x_2+d_3x_4.\end{array}, \leqno (6.10)$$

We deduce from $\delta J_i=\Delta J_i=0$ that $$ \delta_1 J_i=0,
\quad \delta_2 J_i=0, \quad i=1,2.$$

Consider the partial differential equation $\delta J=0$. Solving
it one finds  that $J=K(Y_1,Y_2,Y_3)$ where $K(y_1,y_2,y_3)$ is an
arbitrary function and $$Y_1=x_2, \quad Y_2=x_4, \quad
Y_3=x_4x_1-x_2x_3. \leqno (6.11)$$ Therefore, in view of (6.9),
(6.11) we have $J_1=J_1(Y_1,Y_2,Y_3)$ and $
\displaystyle\frac{\partial J_1}{\partial Y_3}\neq 0.$

Consequently, as $\delta_2 Y_1=\delta_2 Y_2=0$, one gets $$
\delta_2 J_1=\displaystyle \frac{\partial J_1}{\partial
Y_1}\delta_2 Y_1 +\displaystyle \frac{\partial J_1}{\partial
Y_2}\delta_2 Y_2+\displaystyle \frac{\partial J_1}{\partial
Y_3}\delta_2 Y_3=\displaystyle \frac{\partial J_1}{\partial
Y_3}\delta_2 Y_3. $$

This  implies $$ \delta_2 Y_3=0. \leqno (6.12) $$

By substituting in (6.12) the expression  for $Y_3$  given by
(6.11) we arrive to $$ b_3=d_1=0, \quad b_1=d_3=\rho, \leqno
(6.13) $$ for some $\rho \in {\Complex}$.

We now use the equation $\delta_1 J=0$ which can be written as $$
\delta_1 J=\displaystyle \frac{\partial J}{\partial Y_1}\delta_1
Y_1 +\displaystyle \frac{\partial J}{\partial Y_2}\delta_1
Y_2+\displaystyle \frac{\partial J}{\partial Y_3}\delta_1 Y_3=0,
\leqno (6.14) $$

One can show that $$
\begin{array}{lll}
\delta_1Y_1=\rho Y_1,\\ \delta_1Y_2=\rho Y_2, \\
\delta_1Y_3=v_{1}Y_1^2+v_2Y_2^2+v_3Y_1Y_2, \end{array} $$ where
$v_1=d_2-c_1$, $v_2=a_3-b_4$, $v_3=a_1-b_2-c_3+d_4$.

Hence, (6.14) yields $$ \rho Y_1 \displaystyle \frac{\partial
J}{\partial Y_1} +\rho Y_2\displaystyle \frac{\partial J}{\partial
Y_2}+(v_{1}Y_1^2+v_2Y_2^2+v_3Y_1Y_2)\displaystyle \frac{\partial
J}{\partial Y_3}=0.$$

This equation possesses  two rational functionally independent
solutions $J_1(Y_1,Y_2,Y_3)$, $J_2(Y_1,Y_2,Y_3)$ only if $$
\rho=0, \quad v_1=v_2=v_3=0 .$$ which gives $$ a_1=\epsilon_1+b_2,
\quad c_3=\epsilon_1+d_4,\quad c_1=d_2=\zeta_1, \quad
a_3=b_4=\zeta_2, \quad \epsilon_1, \zeta_1, \zeta_2 \in
{\Complex}. \leqno (6.15) $$

After  substitutions  of (6.13), (6.15) in (6.2)  the matrix $R$
is written as $$ R=\left(
\begin{array}{cccc}
b_2+\epsilon_1 & a_2 & \zeta_2 & a_4 \\ 0 & b_2 & 0 & \zeta_2 \\
\zeta_1 & c_2 & d_4+\epsilon_1 & c_4
\\ 0 & \zeta_1 & 0 & d_4
\end{array}
\right).$$

Now, consider the characteristic polynomial $P(\lambda)$ of $R$
$$P(\lambda)=\lambda^4+P_1\lambda^3+P_2\lambda^2+P_3\lambda+P_4,$$
where $$
\begin{array}{llll}
P_1=-2(b_2+d_4+\epsilon_1),\\ P_2=3b_2\epsilon_1-2\zeta_1\zeta_2+
3\epsilon_1d_4+4b_2d_4+b_2^2+\epsilon^2_1 +d_4^2,\\
P_3=-(d_4+b_2+\epsilon_1)(2b_2d_4+b_2\epsilon_1+d_4\epsilon_1-2\zeta_1
\zeta_2), \\ P_4=(b_2d_4
-\zeta_1\zeta_2)(b_2d_4+b_2\epsilon_1+d_4\epsilon_1-\zeta_1\zeta_2+\epsilon_1^2).
\end{array}$$
As above, in view of (6.1), (6.2) all eigenvalues of $R$ must be
equal to $0$ and therefore  $P_i=0$, $1 \leq i \leq 4$. This
system gives $$ \epsilon_1=0, \quad b_2=\eta_1, \quad d_4=-\eta_1,
\quad \eta_1^2+\zeta_1 \zeta_2=0,$$ and the monodromy matrix $T_2$
becomes
$$ T_2=\left(
\begin{array}{cccc}
\eta_1+1 & a_2 & \zeta_2 & a_4 \\ 0 & \eta_1+1 & 0 & \zeta_2 \\
\zeta_1 & c_2 & 1-\eta_1 & c_4
\\ 0 & \zeta_1 & 0 & 1-\eta_1
\end{array}
\right). $$

The matrix  $T_1T_2$ has the eigenvalues $\{1,1,1,1\}$ which
contradicts to c) of Lemma 5.1 and proves our claim. \quad \quad
$\Box$

Due to our definition of integrability we deduce from Theorem 6.2
the following

\vspace{0.5cm}

\noindent {\bf Theorem 6.3} {\it The planar three-body problem is
meromorphically non-integrable near the Lagrangian parabolic
solution. }

\begin{center}

{\bf 7. Final remarks}
\end{center}

In the end of 19th century Poincar\'e [18] indicated some
qualitative phenomena in the behavior of phase trajectories which
prevent the appearance of new integrals of a Hamiltonian system
besides those which are present, but fail to form a set sufficient
for complete integrability.

Let $M^{2n}$ be the phase space, and $H:M^{2n}\rightarrow \Real$,
$H=H_0+\epsilon H_1 +O(\epsilon^2)$ the Hamiltonian function.
Suppose that for $\epsilon =0$ the corresponding Hamiltonian
system has an $m$--dimensional hyperbolic invariant torus $T^m_0$.
According to the Graff's theorem [5], for small $\epsilon$ the
perturbed system has an invariant hyperbolic torus
$T^m_{\epsilon}$ depending analytically on $\epsilon$. It can be
shown that  $T^m_{\epsilon}$ has asymptotic invariant manifolds
$\Lambda^+$ and $\Lambda^-$ filled with trajectories which tend to
the torus $T^m_{\epsilon}$ as $t\rightarrow +\infty$ and
$t\rightarrow -\infty$ respectively. In integrable Hamiltonian
systems such manifolds (called also {\it separatrices}), as a
rule, coincide. In the nonintegrable cases, the situation is
different: asymptotic surfaces can have transverse intersection
forming a complicated tangle which prevent the appearance of new
integrals. For a modern presentation of these results see, for
example, [7].

The method of splitting of asymptotic surfaces was applied to the
three--body problem by many authors. In his book [13] J.K. Moser
described a technique which use the symbolic dynamics associated
with a transverse homoclinic point. Applying this method, it was
shown in [9] that under certain assumptions the planar circular
restricted three--body problem does not possess an additional real
analytic integral. The similar result for the Sitnikov problem and
the collinear three--body problem can be found in [13], [10]. The
existence and the transverse intersection of stable and unstable
manifolds along some periodic orbits in the planar three--body
problem where two masses are sufficiently small was established in
[11], using the  results obtained in [12].

It is necessary to note that Theorem 6.3  implies the nonexistence
of a complete set of complex analytic first integrals for the
general planar three--body problem. To prove the nonexistence of
real analytic integrals one should use some heteroclinic
phenomena and can propose the following line of reasoning : Let
$M^{\infty}$ be the infinity manifold, then the taken Lagrangian
parabolic orbit is biasymptotic to it. This is a weakly hyperbolic
invariant manifold and the reference orbit is a heteroclinic orbit
to different periodic orbits sitting in $M^{\infty}$. The
dynamical interpretation of Theorem 6.3 seems to be the
transversality of the invariant stable and unstable manifolds of
$M^{\infty}$, along this orbit. A combination of passages near
several of these orbits (there is all the family obtained by
rotation) should allow to prove the existence of a heteroclinic
chain. This, in turn, gives rise to an embedding of a suitable
subshift, with lack of predictability, chaos and implies the
nonexistence of real analytic integrals.

\begin{center}

{\bf  Acknowledgements}

\end{center}
I  would like to thank L. Gavrilov and V. Kozlov for useful
discussions  and the advice to study the present problem. Also, I
thank to J.-P. Ramis, J.J. Morales-Ruiz, J.-A. Weil  and D.
Boucher for their attention to the paper. I am very grateful to
the anonymous referee for his useful remarks.

\begin{center}

{\bf  References}

\end{center}

\noindent [1] V.I. Arnold, V.V. Kozlov, A.I. Neishtadt,{\it
Dynamical systems III, } Springer--Verlag, p. 63, (1987).

\noindent [2] H. Bruns, {\it Ueberdie Integrale des vierk$\ddot
o$rper Problems}, Acta Math. 11, p. 25-96, (1887-1888).

\noindent [3] J. Chazy, { \it Sur l'allure du mouvement dans le
probl\`eme des trois  corps quand le temps croit ind\'efiniment},
Ann. Sci. Ecole Norm., 39 , 29-130 (1922).

\noindent [4] V.V. Golubev, {\it Lectures on analytic theory of
differential equations}, Gostekhizdat, Moskow, (1950), (Russian).

\noindent [5] S.M. Graff, {\it On the conservation of hyperbolic
invariant tori for Hamiltonian systems}, J. Differential Equations
15, 1--69, (1974).

\noindent [6] S.V. Kowalewskaya, {\it Sur le probl\`eme de la
rotation d'un corps solide autour d'un point fixe}, Acta. Math.
12,177-232 (1889).

\noindent [7] V.V. Kozlov, {\it Symmetries, Topology, and
Resonances in Hamiltonian mechanics,} Springer-Verlag (1996).

\noindent [8] J.L. Lagrange, {\it Oeuvres}. Vol. 6, 272-292, Paris
(1873).

\noindent [9] J. Llibre, C. Sim\'o, {\it Oscillatory solutions in
the planar restricted three--body problem,}Math. Ann., 248:
153--184, 1980.

\noindent [10] J. Llibre, C. Sim\'o, {\it Some homoclinic
phenomena in the three--body problem,} J. Differential Equations,
37, no. 3, 444--465, 1980.

\noindent [11] R. Martinez, C. Sim\'o, {\it A note on the
existence of heteroclinic orbits in the planar three body
problem,} In Seminar on Dynamical Systems, Euler International
Mathematical Institute, St. Petersburg, 1991, S. Kuksin, V.
Lazutkin and J. P$\ddot o$chel, editors, 129--139, Birkh$\ddot
a$user, 1993.

\noindent [12] R. Martinez, C. Pinyol, {\it Parabolic orbits in
the elliptic restricted three body problem,} J. Differential
Equations, 111, 299--339, (1994).

\noindent [13] J.K. Moser, {\it Stable and random motions in
dynamical systems,} Princeton Univ. Press, Princeton, N.J., 1973.

\noindent [14] R. Moeckel, {\it Chaotic dynamics near triple
collision}, Arch. Rational. Mech. Anal. 107, no. 1, 37-69 (1989).

\noindent [15] J.J. Morales-Ruiz, J.P. Ramis, {\it Galosian
Obstructions to integrability of Hamiltonian Systems}, Preprint
(1998).

\noindent [16] I. Newton, {\it Philosophiae naturalis principia
mathematica}, Imprimatur S. Pepys, Reg. Soc. Praeses, julii 5,
1686, Londini anno MDCLXXXVII.

\noindent [17] P. Painlev\'e, {\it M\'emoire sur les int\'egrales
premi\'eres du probl\`eme des n corps}, Acta Math. Bull. Astr. T
15 (1898).

\noindent [18] H. Poincar\'e, {\it Les m\'ethodes novelles de la
m\'ecanique c\'eleste}, vol. 1-3. Gauthier--Villars, Paris 1892,
1893, 1899.

\noindent [19] R. Schaefke, Private communication.

\noindent [20] M. Singer, A. Baider, R. Churchill, D. Rod, {\it On
the infinitesimal Geometry of Integrable Systems}, in Mechanics
Day, Shadwich et. al., eds, Fields Institute Communications, 7,
AMS, 5-56 (1996).

\noindent [21] K. F. Sundman, { \it Memoire sur le probl\`eme des
trois corps}, Acta Math. 36, 105-107 (1913).

\noindent [22] C.L. Siegel, J.K. Moser, {\it Lectures on Celestial
Mechanics,} Springer-Verlag (1971).

\noindent [23] Gh. Topan, {\it Sur une int\'egrale premi\'ere
transcendante dans certaines configurations du probl\'eme des
trois corps}, Bull. Math. Soc. Sci. Math. R. S. Roumanie (N.S.),
no. 1, 83-91 (1989).

\noindent [24] E.T. Whittaker, {\it A Treatise on the Analytical
Dynamics of particles and Rigid Bodies}. Cambridge University
Press, New York, (1970).

\noindent [25] A. Winter, {\it The analytical Foundations of
Celestial Mechanics,} Princeton Univ. Press, Princeton, (1941).

\noindent [26] H. Yoshida, {\it A criterion for the nonexistence
of an additional integral in Hamiltonian systems with a
homogeneous potential}, Phys. D. 29, no. 1-2, 128-142, (1987).

\noindent [27] S.L. Ziglin, {\it Branching of solutions and
non-existence of first integrals in Hamiltonian Mechanics I},
Func. Anal. Appl. 16 (1982).

\vspace{0.3cm}

\noindent( please use this address for correspondence )

\noindent{\small \bf
 Section de Mathematiques,\\
 Universit\'e de Gen\`eve\\
 2-4, rue du Lievre,\\
 CH-1211, Case postale 240, Suisse \\
Tel l.: +41 22 309 14 03 \\
 Fax: +41 22 309 14 09 \\
E--mail: Alexei.Tsygvintsev@math.unige.ch}

\vspace{0.3cm}

 \noindent{\small \bf
 Laboratoire Emile Picard, UMR 5580,\\
 Universit\'e Paul Sabatier\\
 118, route de Narbonne,\\
 31062 Toulouse Cedex, France \\
Tel l.: 05 61 55 83 37 \\
 Fax: 05 61 55 82 00 \\
E--mail: tsygvin@picard.ups-tlse.fr}

\vspace{0.3cm}

\noindent {\bf  Tsygvintsev Alexei}

\newpage
\begin{center}
{\bf   Appendix A. The matrices $A_{\infty}$, $A$, $B$, $R$, $J$}
\end{center}

\begin{center}
$A_{\infty}=A_{\infty, ij}$, $1\leq i,j \leq 4$. \\
\end{center}

 \noindent ${A_{\infty , \,11}}={\displaystyle \frac {1}{4}} \,
{\displaystyle \frac {12\,\alpha  + 5\,\beta  + 5\,\beta \,\alpha
 ^{2} + 26\,\alpha \,\beta  + 12\,\alpha ^{2}}{\alpha \,{S_{1}}}
} , \quad {A_{\infty , \,12}}={\displaystyle \frac {3}{4}} \,
{\displaystyle \frac {\sqrt{3}\,(\alpha  + 1)\,\beta \,(\alpha
 - 1)}{\alpha \,{S_{1}}}},$\\
$
{A_{\infty , \,13}}= - 2\,{\displaystyle \frac {{S_{1}}}{{S_{
2}}^{2}\,{S_{3}}^{3}\,\alpha ^{4}\,\beta }}, \quad {A_{\infty ,
\,14}}=0, $\\
$
{A_{\infty , \,21}}={\displaystyle \frac {3}{4}} \, {\displaystyle
\frac {\sqrt{3}\,(\alpha  + 1)\,\beta \,(\alpha
 - 1)}{\alpha \,{S_{1}}}},  \quad {A_{\infty , \,22}}= - {\displaystyle \frac {1}{4}} \,
{\displaystyle \frac { - 12\,\alpha  + \beta  + \beta \,\alpha ^{
2} - 2\,\alpha \,\beta  - 12\,\alpha ^{2}}{\alpha \,{S_{1}}}} ,
$\\
$
{A_{\infty , \,23}}=0, \quad {A_{\infty , \,24}}= -
2\,{\displaystyle \frac {{S_{1}}}{{S_{
2}}^{2}\,{S_{3}}^{3}\,\alpha ^{4}\,\beta }} , $\\
$
{A_{\infty , \,31}}={\displaystyle \frac {1}{8}} \, {\displaystyle
\frac {\alpha ^{2}\,\beta \,{S_{3}}^{3}\,(\alpha
 + 1)\,{S_{2}}^{3}\,(2\,\alpha  + 13\,\beta  + 13\,\beta \,\alpha
 ^{2} + 24\,\alpha \,\beta  + 2\,\alpha ^{2})}{{S_{1}}^{3}}} , $\\
$
{A_{\infty , \,32}}={\displaystyle \frac {3}{8}} \, {\displaystyle
\frac {\sqrt{3}\,(\beta  + 2\,\alpha  + 4\,\alpha \,\beta  + \beta
\,\alpha ^{2} + 2\,\alpha ^{2})\,(\alpha  - 1)\, \beta \,\alpha
^{2}\,{S_{3}}^{3}\,{S_{2}}^{3}}{{S_{1}}^{3}}} , $\\
$
{A_{\infty , \,33}}= - {\displaystyle \frac {1}{4}} \,
{\displaystyle \frac {\beta \,(5\,\alpha ^{2} + 14\,\alpha  + 5)
}{\alpha \,{S_{1}}}} \quad {A_{\infty , \,34}}= - {\displaystyle
\frac {3}{4}} \, {\displaystyle \frac {\sqrt{3}\,(\alpha  +
1)\,\beta \,(\alpha
 - 1)}{\alpha \,{S_{1}}}} , $\\
$
{A_{\infty , \,41}}={\displaystyle \frac {3}{8}} \, {\displaystyle
\frac {\sqrt{3}\,(\beta  + 2\,\alpha  + 4\,\alpha \,\beta  + \beta
\,\alpha ^{2} + 2\,\alpha ^{2})\,(\alpha  - 1)\, \beta \,\alpha
^{2}\,{S_{3}}^{3}\,{S_{2}}^{3}}{{S_{1}}^{3}}} , $\\
$
{A_{\infty , \,42}}={\displaystyle \frac {1}{8}} \, {\displaystyle
\frac {\alpha ^{2}\,\beta \,{S_{3}}^{3}\,(\alpha
 + 1)\,{S_{2}}^{3}\,( - 10\,\alpha  + 7\,\beta  + 7\,\beta \,
\alpha ^{2} - 12\,\alpha \,\beta  - 10\,\alpha ^{2})}{{S_{1}}^{3}
}}, $\\
$
{A_{\infty , \,43}}= - {\displaystyle \frac {3}{4}} \,
{\displaystyle \frac {\sqrt{3}\,(\alpha  + 1)\,\beta \,(\alpha
 - 1)}{\alpha \,{S_{1}}}} , \quad {A_{\infty , \,44}}={\displaystyle \frac {1}{4}} \,
{\displaystyle \frac {\beta \,(10\,\alpha  + \alpha ^{2} + 1)}{
\alpha \,{S_{1}}}}. $

\begin{center}
 $A=(A_{ij})$, $1\leq i,j \leq 4$.\\
\end{center}

\noindent$ {A_{11}}= - {\displaystyle \frac {1}{4}}
\,{\displaystyle \frac {(\alpha  + 1)\,(\alpha \,\beta  +
4\,\alpha  + \beta )}{ \alpha \,{S_{1}}}} , \quad
{A_{12}}={\displaystyle \frac {1}{4}} \,{\displaystyle \frac
{\sqrt{3}\,(\alpha  + 1)\,\beta \,(\alpha  - 1)}{\alpha \,{
S_{1}}}} , $\\ $ {A_{13}}=2\,{\displaystyle \frac
{{S_{1}}}{{S_{2}}^{2}\,{S_{3 }}^{3}\,\alpha ^{4}\,\beta }},\quad
{A_{14}}=0,\quad {A_{21}}={\displaystyle \frac {1}{4}}
\,{\displaystyle \frac {\sqrt{3}\,(\alpha  + 1)\,\beta \,(\alpha -
1)}{\alpha \,{ S_{1}}}} , $\\
$
{A_{22}}= - {\displaystyle \frac {1}{4}} \,{\displaystyle \frac
{10\,\alpha \,\beta  + 3\,\beta \,\alpha ^{2} + 3\,\beta
 + 4\,\alpha ^{2} + 4\,\alpha }{\alpha \,{S_{1}}}} ,\quad {A_{24}}=2\,{\displaystyle \frac {{S_{1}}}{{S_{2}}^{2}\,{S_{3
}}^{3}\,\alpha ^{4}\,\beta }}$\\
$
{A_{31}}= - {\displaystyle \frac {1}{8}} \,{\displaystyle \frac
{\alpha ^{2}\,\beta ^{2}\,{S_{3}}^{3}\,(\alpha  + 1)\,( \alpha  -
1)^{2}\,{S_{2}}^{3}}{{S_{1}}^{3}}} ,\quad {A_{34}}= -
{\displaystyle \frac {1}{4}} \,{\displaystyle \frac
{\sqrt{3}\,(\alpha  + 1)\,\beta \,(\alpha  - 1)}{\alpha \,{
S_{1}}}},$\\
$
A_{23}=0, \quad {A_{32}}={\displaystyle \frac {1}{8}}
\,{\displaystyle \frac {\sqrt{3}\,(\alpha  - 1)\,(\alpha  +
1)^{2}\,\alpha ^{2}\, \beta
^{2}\,{S_{2}}^{3}\,{S_{3}}^{3}}{{S_{1}}^{3}}} ,\quad
{A_{33}}={\displaystyle \frac {1}{4}} \,{\displaystyle \frac
{(\alpha - 1)^{2}\,\beta }{\alpha \,{S_{1}}}},$\\
${A_{41}}={\displaystyle \frac {1}{8}} \,{\displaystyle \frac
{\sqrt{3}\,(\alpha  - 1)\,(\alpha  + 1)^{2}\,\alpha ^{2}\, \beta
^{2}\,{S_{2}}^{3}\,{S_{3}}^{3}}{{S_{1}}^{3}}} ,\quad {A_{42}}= -
{\displaystyle \frac {3}{8}} \,{\displaystyle \frac {\alpha
^{2}\,\beta ^{2}\,{S_{3}}^{3}\,(\alpha  + 1)^{3}\,{
S_{2}}^{3}}{{S_{1}}^{3}}} , $\\
$
{A_{43}}= - {\displaystyle \frac {1}{4}} \,{\displaystyle \frac
{\sqrt{3}\,(\alpha  + 1)\,\beta \,(\alpha  - 1)}{\alpha \,{
S_{1}}}} ,\quad {A_{44}}={\displaystyle \frac {3}{4}}
\,{\displaystyle \frac {(\alpha  + 1)^{2}\,\beta }{\alpha
\,{S_{1}}}}. $

\newpage

\begin{center}
$R=(R_{ij})$, $1\leq i,j \leq 4$.\\
\end{center}

\noindent $ {R_{11}}= - {\displaystyle \frac {1}{2}}
\,{\displaystyle \frac {2\,\alpha  + \beta  + 6\,\alpha \,\beta  +
2\,\alpha ^{2}
 + \beta \,\alpha ^{2}}{\alpha \,{S_{1}}}}, \quad {R_{12}}= - {\displaystyle \frac {1}{2}} \,{\displaystyle
\frac {\sqrt{3}\,(\alpha  + 1)\,\beta \,(\alpha  - 1)}{\alpha \,{
S_{1}}}} , $\\
$
{R_{13}}=0, \quad {R_{14}}=0, $\\
$
{R_{21}}= - {\displaystyle \frac {1}{2}} \,{\displaystyle \frac
{\sqrt{3}\,(\alpha  + 1)\,\beta \,(\alpha  - 1)}{\alpha \,{
S_{1}}}},\quad {R_{22}}={\displaystyle \frac {1}{2}}
\,{\displaystyle \frac {(\alpha  + 1)\,( - 2\,\alpha  + \alpha
\,\beta  + \beta ) }{\alpha \,{S_{1}}}} , $\\
$
{R_{23}}=0, \quad {R_{24}}=0, $\\
$
{R_{31}}= - {\displaystyle \frac {1}{8}} \,{\displaystyle \frac
{\beta \,\alpha ^{2}\,{S_{2}}^{3}\,{S_{3}}^{3}\,(\alpha  +
1)\,(\alpha ^{2} + 6\,\beta \,\alpha ^{2} + \alpha  + 13\,\alpha
\,\beta  + 6\,\beta )}{{S_{1}}^{3}}} , $\\
$
{R_{32}}= - {\displaystyle \frac {1}{8}} \,{\displaystyle \frac
{\sqrt{3}\,(3\,\alpha ^{2} + 2\,\beta \,\alpha ^{2} + 7\, \alpha
\,\beta  + 3\,\alpha  + 2\,\beta )\,(\alpha  - 1)\,\beta \,\alpha
^{2}\,{S_{3}}^{3}\,{S_{2}}^{3}}{{S_{1}}^{3}}} , $\\
$
{R_{33}}={\displaystyle \frac {1}{2}} \,{\displaystyle \frac
{(\alpha  + \sqrt{3} + 2)\,(\alpha  + 2 - \sqrt{3})\,\beta
}{\alpha \,{S_{1}}}} ,\quad {R_{34}}={\displaystyle \frac {1}{2}}
\,{\displaystyle \frac {\sqrt{3}\,(\alpha  + 1)\,\beta \,(\alpha -
1)}{\alpha \,{ S_{1}}}}, $\\
$
{R_{41}}= - {\displaystyle \frac {1}{8}} \,{\displaystyle \frac
{\sqrt{3}\,(3\,\alpha ^{2} + 2\,\beta \,\alpha ^{2} + 7\, \alpha
\,\beta  + 3\,\alpha  + 2\,\beta )\,(\alpha  - 1)\,\beta \,\alpha
^{2}\,{S_{3}}^{3}\,{S_{2}}^{3}}{{S_{1}}^{3}}}, $\\
$
{R_{42}}= - {\displaystyle \frac {1}{8}} \,{\displaystyle \frac
{\beta \,\alpha ^{2}\,{S_{2}}^{3}\,{S_{3}}^{3}\,(\alpha  + 1)\,( -
5\,\alpha ^{2} + 2\,\beta \,\alpha ^{2} - 5\,\alpha  - 9 \,\alpha
\,\beta  + 2\,\beta )}{{S_{1}}^{3}}} , $\\
$
{R_{43}}={\displaystyle \frac {1}{2}} \,{\displaystyle \frac
{\sqrt{3}\,(\alpha  + 1)\,\beta \,(\alpha  - 1)}{\alpha \,{
S_{1}}}}, \quad {R_{44}}= - {\displaystyle \frac {1}{2}}
\,{\displaystyle \frac {(\alpha  + \sqrt{3} + 2)\,(\alpha  + 2 -
\sqrt{3})\,\beta }{\alpha \,{S_{1}}}}. $\\

\begin{center}
$J=(J_{ij})$, $1\leq i,j \leq 4$.\\
\end{center}

\noindent $ {J_{11}}= - {\displaystyle \frac {1}{2}}
\,{\displaystyle \frac {\sqrt{3}\,(\alpha  + 1)\,\beta \,(\alpha -
1)}{\alpha \,{ S_{1}}}} , \quad {J_{12}}={\displaystyle \frac
{1}{2}} \,{\displaystyle \frac {(\alpha  + 1)\,( - 2\,\alpha  +
\alpha \,\beta  + \beta ) }{\alpha \,{S_{1}}}}, $\\
$
{J_{13}}=0 ,\quad {J_{14}}=0, $\\
$
{J_{21}}={\displaystyle \frac {1}{2}} \,{\displaystyle \frac
{2\,\alpha  + \beta  + 6\,\alpha \,\beta  + 2\,\alpha ^{2}
 + \beta \,\alpha ^{2}}{\alpha \,{S_{1}}}}, \quad {J_{22}}={\displaystyle \frac {1}{2}} \,{\displaystyle
\frac {\sqrt{3}\,(\alpha  + 1)\,\beta \,(\alpha  - 1)}{\alpha \,{
S_{1}}}}, $\\
$
{J_{23}}=0 , \quad {J_{24}}=0, \quad {J_{31}}= - {\displaystyle
\frac {1}{4}} \,{\displaystyle \frac {\sqrt{3}\,(\alpha  -
1)\,(\alpha  + 1)^{2}\,\alpha ^{2}\, \beta
^{2}\,{S_{2}}^{3}\,{S_{3}}^{3}}{{S_{1}}^{3}}} , $\\
$
{J_{32}}={\displaystyle \frac {1}{4}} \,{\displaystyle \frac
{\beta ^{2}\,\alpha ^{2}\,{S_{2}}^{3}\,{S_{3}}^{3}\,(\alpha
  + 1)\,(\alpha ^{2} + 4\,\alpha  + 1)}{{S_{1}}^{3}}} , $\\
$
{J_{33}}={\displaystyle \frac {1}{2}} \,{\displaystyle \frac
{\sqrt{3}\,(\alpha  + 1)\,\beta \,(\alpha  - 1)}{\alpha \,{
S_{1}}}} , \quad {J_{34}}= - {\displaystyle \frac {1}{2}}
\,{\displaystyle \frac {2\,\alpha  + \beta  + 6\,\alpha \,\beta  +
2\,\alpha ^{2}
 + \beta \,\alpha ^{2}}{\alpha \,{S_{1}}}}, $\\
$
{J_{41}}={\displaystyle \frac {1}{4}} \,{\displaystyle \frac
{\beta ^{2}\,\alpha ^{2}\,{S_{2}}^{3}\,{S_{3}}^{3}\,(\alpha
  + 1)\,(\alpha ^{2} + 4\,\alpha  + 1)}{{S_{1}}^{3}}} , $\\
$
{J_{42}}={\displaystyle \frac {1}{4}} \,{\displaystyle \frac
{\sqrt{3}\,(\alpha  - 1)\,(\alpha  + 1)^{2}\,\alpha ^{2}\, \beta
^{2}\,{S_{2}}^{3}\,{S_{3}}^{3}}{{S_{1}}^{3}}} , $\\
$
{J_{43}}= - {\displaystyle \frac {1}{2}} \,{\displaystyle \frac
{(\alpha  + 1)\,( - 2\,\alpha  + \alpha \,\beta  + \beta )
}{\alpha \,{S_{1}}}}, \quad {J_{44}}= - {\displaystyle \frac
{1}{2}} \,{\displaystyle \frac {\sqrt{3}\,(\alpha  + 1)\,\beta
\,(\alpha  - 1)}{\alpha \,{ S_{1}}}}. $

\end{document}